\documentclass{article}

\usepackage{arxiv}

\usepackage[utf8]{inputenc} 
\usepackage[T1]{fontenc}    
\usepackage{hyperref}       
\usepackage{url}            
\usepackage{booktabs}       
\usepackage{amsfonts}       
\usepackage{nicefrac}       
\usepackage{microtype}      
\usepackage{lipsum}		
\usepackage{graphicx}
\usepackage{natbib}
\setcitestyle{numbers,square}
\usepackage{doi}
\usepackage{amsmath}
\usepackage{subcaption}
\usepackage{authblk}
\usepackage{orcidlink}

\title{A CONVEXIFICATION-BASED OUTER-APPROXIMATION METHOD FOR CONVEX AND NONCONVEX MINLP}

\author[1]{Zedong Peng \orcidlink{0000-0001-6001-1738}}
\author[1]{Kaiyu Cao  \orcidlink{0000-0002-8901-6310}}
\author[2]{Kevin C. Furman}
\author[1]{Can Li}
\author[3]{Ignacio E. Grossmann \orcidlink{0000-0002-7210-084X}}
\author[1,4,5]{David E. Bernal Neira \orcidlink{0000-0002-8308-5016}}

\affil[1]{Davidson School of Chemical Engineering, Purdue University, 480 Stadium Mall Drive, West Lafayette, IN, 47907, USA}
\affil[2]{ExxonMobil Upstream Research Co, 22777 Springwoods Village Parkway, Spring, TX, 77389, USA}
\affil[3]{Department of Chemical Engineering, Carnegie Mellon University, Pittsburgh, PA 15213, USA}
\affil[4]{Research Institute for Advanced Computer Science, Universities Space Research Association, Mountain View, CA, USA}
\affil[5]{Quantum Artificial Intelligence Laboratory (QuAIL), NASA Ames Research Center, Mountain View, CA, USA}

\renewcommand{\shorttitle}{A Convexification-based Outer-Approximation Method for Convex and Nonconvex MINLP}

\hypersetup{
pdftitle={A CONVEXIFICATION-BASED OUTER-APPROXIMATION METHOD FOR CONVEX AND NONCONVEX MINLP},
}

\begin{document}
\maketitle

\begin{abstract}
The advancement of domain reduction techniques has significantly enhanced the performance of solvers in mathematical programming. This paper delves into the impact of integrating convexification and domain reduction techniques within the Outer- Approximation method. We propose a refined convexification-based Outer-Approximation method alongside a Branch-and-Bound method for both convex and nonconvex Mixed-Integer Nonlinear Programming problems. These methods have been developed and incorporated into the open-source Mixed-Integer Nonlinear Decomposition Toolbox for Pyomo-MindtPy. Comprehensive benchmark tests were conducted, validating the effectiveness and reliability of our proposed algorithms. These tests highlight the improvements achieved by incorporating convexification and domain reduction techniques into the Outer-Approximation and Branch-and-Bound methods.
\end{abstract}


\section{Introduction}
Mixed-integer nonlinear programming (MINLP) has broad applications in process systems engineering (PSE), including planning, scheduling, and control. It offers a powerful modeling framework that optimizes discrete and continuous variables involved in linear and nonlinear constraints. However, the combinatorial complexity, nonlinearity, and even nonconvexity lead to substantial challenges in optimizing such problems. 

Generally, MINLP can be classified as convex and nonconvex, depending on the convexity of its continuous relaxation. The algorithms for MINLP are primarily categorized into Branch-and-Bound (B\&B) methods and decomposition methods \cite{kronqvist2019review}. The main idea of MINLP decomposition algorithms is to generate linear inequalities to approximate nonlinear constraints and iteratively solve the relaxed Mixed-Integer Linear Programming (MILP) main problem and Nonlinear Programming (NLP) subproblems. Decomposition methods for convex MINLP problems include the Outer-Approximation (OA) method \cite{duran1986outer}. This method involves solving an MILP defined by linear inequalities that relax the nonlinear constraints, known as OA cuts, and an NLP with the main problem's integer solution fixed. The first-order Taylor approximation of the nonlinear constraints defines these OA cuts. To reduce the MILP problem solution time, the Linear Programming and Nonlinear Programming-based B\&B (LP/NLP-B\&B) method \cite{quesada1992lp} maintains the same B\&B for the MILP, solves NLPs at the tree's integer nodes, and uses OA cuts to improve the searching bounds. Consequently, these decomposition methods are often called multi-tree and single-tree, based on their management of the MILP problem. A significant limitation of these methods is their initialization, usually given by relaxing the nonlinearity of the problem and relying on the cuts generated at each iteration to provide a better linear approximation as iterations progress.

While the OA method is effective for convex MINLPs, its limitations become apparent for nonconvex MINLP, where OA cuts do not guarantee validity in relaxations of nonlinear functions, precluding global optimality guarantees. Other relaxation techniques have been developed for nonconvex MINLPs addressing this challenge. Among these, the Auxiliary Variable Method (AVM) and McCormick relaxations are successful strategies for generating relaxations of nonconvex factorable functions \cite{tawarmalani2013convexification}. AVM achieves this by introducing an auxiliary variable and a corresponding equality constraint for each intermediate nonlinear factor in a function, leading to computational efficiency by decomposing the function into simpler, lower-dimensional components. However, this method entails incorporating many auxiliary variables and constraints. In contrast, McCormick relaxations maintain the dimension of the original function and use a recursive approach to produce the required convex and concave relaxations effectively \cite{Chachuat2013mc}.

In addition to advancements in reformulations and optimization algorithms, the performance of optimization solvers has significantly improved through domain reduction techniques. These techniques encompass bound tightening, eliminating redundant variables and constraints, and convexification \cite{zhang2020optimality}. The bound tightening techniques include Feasibility-based Bound Tightening (FBBT), Optimality-based Bound Tightening (OBBT), and Marginals-based Bound Tightening \cite{zhang2020optimality}. Domain reduction methods, including convexification cuts and bound tightening techniques, have been successful for spatial B\&B methods. These tighter relaxations provide stronger dual bounds, accelerating the B\&B process by facilitating node pruning and efficiently identifying optimal solutions. However, decomposition-based MINLP solvers have not yet fully harnessed the potential of domain reduction techniques.

This work investigates the efficacy of domain reduction techniques within OA methods, applicable to both convex and nonconvex MINLP problems. It is understood that domain reduction techniques can be used during the presolve stage and at each node within the B\&B tree, a strategy known as the branch-and-reduce method. Similarly, in the OA method, these techniques can be employed both in the presolve phase and during solving integer-fixed NLP subproblems. However, this work focuses on the impact of domain reduction methods at the method's initialization stage.

\section{Solution algorithm}

The general form of a MINLP problem is as follows.

\begin{equation}
\label{eq:MINLP}
\tag{MINLP}
\begin{split}
   \min_{x,y} \quad &f(x,y) \\
    s.t. \quad & g_j\left(x,y\right) \le 0,\forall j=1,\ldots,l \\
    & x \in \left[\underline{x}, \overline{x} \right] \subseteq \mathbb{R}^n, \\
    & y \in \left\{\underline{y}, \ldots, \overline{y} \right\} \subseteq \mathbb{Z}^m, \\ 
\end{split}
\end{equation}
where $x$ and $y$ represent continuous variables and discrete variables, respectively. Upper and lower variable bounds are determined by over- and underbars, respectively. Both the objective $f(x,y)$ and constraints $g_j (x,y)$ are potentially nonlinear functions.
The OA solution method for MINLP involves an iterative two-step procedure. The first step in iteration k determines the integer variables' values $y^{k+1}$ by solving problem \eqref{eq:OA-master}, defined by the OA cuts. Its optimal objective function value, encoded in variable $\mu$, provides a dual lower bound (LB) to the original MINLP problem's optimal objective.

\begin{equation}
\label{eq:OA-master}
\tag{OA-MILP}
\begin{split}
    \min_{x, y, \mu} \quad & \mu \\
    s.t. \quad & f\left(x^i,y^i\right)+\nabla f\left(x^i,y^i\right) ^\mathrm{T}
    \left[ \begin{matrix} x - x^i\\ y - y^i \end{matrix} \right]
    \le \mu \quad \forall i = 1, \ldots, k \\
    \quad & g_j\left(x^i,y^i\right)+\nabla g_j\left(x^i,y^i\right) ^\mathrm{T}
    \left[ \begin{matrix} x - x^i\\ y - y^i \end{matrix} \right]
    \le 0 \quad \forall i = 1, \ldots, k, \forall j \in I_i \\
    & x \in \left[\underline{x}, \overline{x} \right] \subseteq \mathbb{R}^n, \\
    & y \in \left\{\underline{y}, \ldots, \overline{y} \right\} \subseteq \mathbb{Z}^m, \\ 
    & \mu \in \mathbb{R}^1.
\end{split}
\end{equation}

The second step is determining the continuous variables' values $x^{k+1}$ by solving problem \eqref{eq:NLP-I} whose optimal solution yields a primal upper bound (UB) to problem (MINLP).

\begin{equation}
\label{eq:NLP-I}
\tag{NLP-I}
\begin{split}
   \min_{x} \quad &f(x,y^{k+1}) \\
    s.t. \quad & g_j\left(x,y^{k+1}\right) \le 0 \quad \forall j=1,\ldots,l \\
    & x \in \left[\underline{x}, \overline{x} \right] \subseteq \mathbb{R}^n. \\
\end{split}
\end{equation}

If problem \eqref{eq:NLP-I} is infeasible, the following feasibility subproblem is solved to minimize a norm $p$ of the constraint violations $s$, as a result of updating $x^{k+1}$.

\begin{equation}
\label{eq:NLP-f}
\tag{NLP-f}
\begin{split}
   \min_{x,s} \quad & \| s\|_{p} \\
    s.t. \quad & g_j\left(x,y^{k+1}\right) \le s_j \quad \forall j=1,\ldots,l \\
    & x \in \left[\underline{x}, \overline{x} \right] \subseteq \mathbb{R}^n. \\
\end{split}
\end{equation}

\begin{figure}[htbp]
\includegraphics[width=0.55\linewidth]{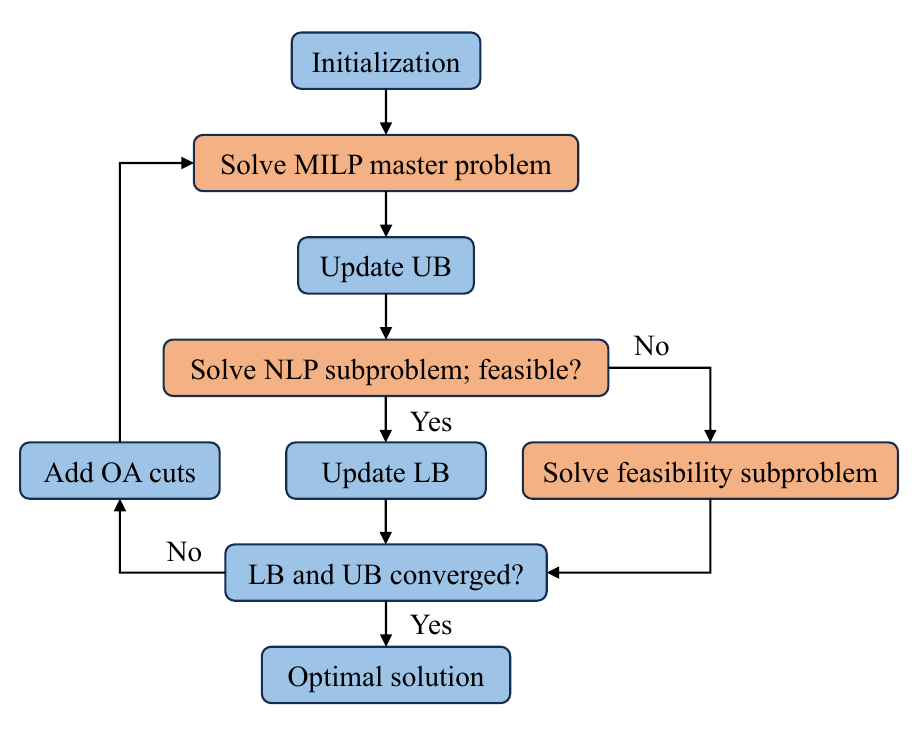}
\centering
\caption{Outer-Approximation method}
\label{fig:OA}
\end{figure}

As shown in Figure \ref{fig:OA}, the OA method begins by solving the relaxed NLP problem and then iteratively solves the \eqref{eq:OA-master}, \eqref{eq:NLP-I}, and \eqref{eq:NLP-f} problems. The key to this process is the progressive accumulation of OA cuts, which incrementally narrows the gap between the LB and the UB. The iterations continue until LB and UB converge, culminating in the OA method reaching the optimal solution. This method is guaranteed to find the global optimal solution of convex MINLPs \cite{duran1986outer}.

\begin{figure}[htbp]
\includegraphics[width=0.8\linewidth]{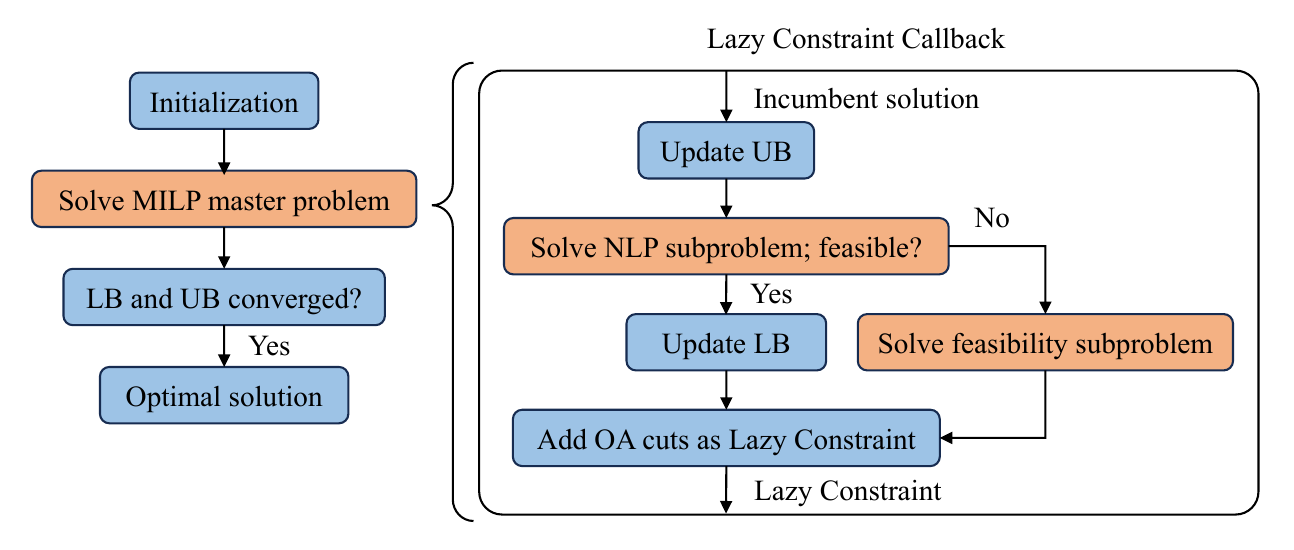}
\centering
\caption{LP/NLP-based Branch and Bound method}
\label{fig:LP-NLP BB}
\end{figure}

Maintaining a single MILP tree for the LP/NLP-B\&B method can be implemented using the LazyConstraint callback function through callback functions in current MILP solvers, as shown in Figure \ref{fig:LP-NLP BB}. This method initializes by solving the relaxed NLP problem as well. Then, a B\&B method is used to solve problem \eqref{eq:OA-master}. Whenever an incumbent solution is found in the search tree, \eqref{eq:NLP-I} is solved, and OA cuts are added as lazy constraints to the MILP tree. This B\&B process is guaranteed to terminate at the global optimal solution of convex MINLP problems \cite{quesada1992lp}.

The OA and LP/NLP-B\&B methods generate tight cuts at the boundary of the nonlinear feasible region defined by the original problem constraints by incurring the cost of solving NLP subproblems. If the \eqref{eq:NLP-I} subproblem is infeasible, \eqref{eq:NLP-f} is solved to find the point closest to the feasible region to generate the tightest possible cuts. Since the linear inequality constraints are accumulated iteratively, the main problem in the early iterations is that a poor approximation of the original MINLP model is obtained. Consequently, the integer combination provided by the main problem tends to yield an infeasible \eqref{eq:NLP-I} subproblem, and no primal bound can be obtained.

In this work, we apply the domain reduction methods to the initialization stage of the OA method to resolve this issue. Eq. \eqref{eq:tightened continuous variables} and \eqref{eq:tightened integer variables} shows the tightened bounds of the discrete and continuous variables, $\left[\underline{x}^\prime,\overline{x}^\prime \right]$ and $\{\underline{y}^\prime,\ldots,\overline{y}^\prime \}$, respectively. Eq. \eqref{eq:convexification cuts} corresponds to the convexification linear cuts generated by the auxiliary variable method or reformulations, where z are auxiliary variables. Since all Eq. \eqref{eq:tightened continuous variables} - \eqref{eq:convexification cuts} are applied at the initialization stage, they are valid for \eqref{eq:OA-master}, \eqref{eq:NLP-I}, and \eqref{eq:NLP-f} problems.

\begin{equation}
    \label{eq:tightened continuous variables}
    x\in \left[{\underline{x}}^\prime,\overline{x}^\prime\right]\in\left[\underline{x},\overline{x}\right]\subseteq\mathbb{R}^n
\end{equation}

\begin{equation}
    \label{eq:tightened integer variables}
    y\in\left\{\underline{y}^\prime,\ldots,\overline{y}^\prime\right\}\in\left\{\underline{y},\ldots,\overline{y}\right\}\subseteq\mathbb{Z}^m
\end{equation}

\begin{equation}
    \label{eq:convexification cuts}
    Ax + By + Cz \le b
\end{equation}

Considering that the convexification cuts are the relaxation of nonlinear constraints, they are redundant in the NLP subproblems where the original nonlinear constraints are included. Therefore, we denote the \eqref{eq:NLP-I} and \eqref{eq:NLP-f} problems with convexification cuts and tightened bounds as complete-scale NLP problems. The \eqref{eq:NLP-I} and \eqref{eq:NLP-f} problems with only tightened bounds are denoted reduced-scale NLP problems.

This work also considers modified alternatives of the OA and LP/NLP-B\&B methods to guarantee global optimality for nonconvex MINLP problems. Several modifications have been introduced to provide such global optimality guarantees, denoted as global OA (GOA) and global LP/NLP-B\&B (GLP/NLP-B\&B). First, instead of adding OA cuts, the affine underestimators and overstimators are generated based on the convex and concave McCormick relaxations using subgradient propagation \cite{Chachuat2013mc}. Second, to guarantee the algorithm's convergence, no-good cuts are generated to cut off the explored integer combinations and prevent the algorithm from repeatedly cycling through the same combinations. These enhancements enable the global algorithms to converge to the global optimum of nonconvex MINLP problems if the NLP subproblems are solved to global optimality \cite{kesavan2004outer}. Furthermore, we integrate the domain reduction techniques in GOA and GLP/NLP-B\&B and investigate their effect on its performance.

\section{Benchmarking and results}
To evaluate the impact of domain reduction techniques, we use test instances from the problem library MINLPLib \cite{vigerske2014towards}. 434 convex instances and 181 nonconvex instances are selected, adhering to the criteria that each instance must have at least one discrete variable and at least one continuous variable. For clarity in our analysis, we use (r) and (c) to distinguish between the reduced-scale and complete-scale NLP subproblems used in the convexification-based OA and LP/NLP-B\&B methods. Moreover, we indicate the results with convexification with the prefix (C-) in the following results.

The benchmark implementation is based on the Mixed-integer nonlinear decomposition toolbox for Pyomo-MindtPy \cite{bernal2018mixed}. We use both the multi-tree and single-tree implementation of OA and GOA strategy, maintaining their default configurations as a baseline. Moreover, a special version of BARON 19.4.4 is used to tighten the bounds and generate convexification cuts. Option dolocal is set to 0, and numloc is set to 0 to turn off local search during upper bounding and preprocessing in BARON. All range reduction and relaxation options are retained at their default settings. Nonlinear FBBT, OBBT, marginals-based, and linear-feasibility-based bound tightening are applied. Outer approximations of convex univariate functions and cutting planes are also applied. For the termination criteria of the algorithm, we set the absolute tolerances $\epsilon_{abs}=10^{-5}$ and relative tolerances $\epsilon_{rel}=10^{-3}$, along with a time limit of 900s. We use GUROBI 10.0.0 as the MILP solver, IPOPTH 3.14 as the NLP solver for convex instances, and BARON 23.6.22 as the NLP solver for nonconvex instances. All tests ran on a Linux cluster with 48 AMD EPYC 7643 2.3GHz CPUs and 1 TB RAM, with each test restricted to using only a single thread.

\begin{figure}
    \centering
    \includegraphics[width=\linewidth]{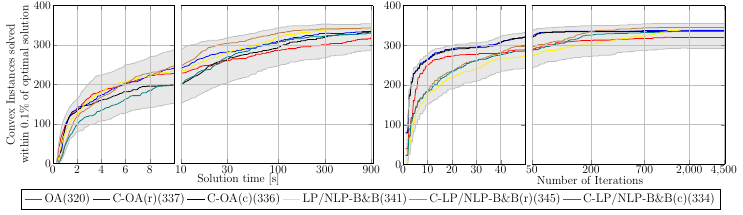}
    \caption{Time and iteration performance profile of 434 convex MINLP instances}
    \label{fig:convex performance profiles}
\end{figure}

The time and iteration performance profiles of the convex instances are presented in Figure \ref{fig:convex performance profiles}. For the LP/NLP-B\&B method, the number of iterations refers to the number of \eqref{eq:NLP-I} subproblems solved. Overall, the convexification-based OA and LP/NLP-B\&B methods utilizing reduced-scale NLP subproblems outperform the other solver alternatives regarding solution time. Regarding the number of iterations, both the OA and LP/NLP-B\&B methods benefit from the convexification cuts and the bound tightening techniques. Interestingly, the choice between complete-scale and reduced-scale NLP subproblems does not significantly impact iteration performance.

However, it is noteworthy that the convexification-based LP/NLP-B\&B method with complete-scale NLP subproblems underperforms in time performance compared to the standard LP/NLP-B\&B method. This coincides with our previous statement that the convexification cuts are redundant and increase the computational complexity of \eqref{eq:NLP-I} and \eqref{eq:NLP-f} subproblems. For instances that can be solved within one second, both the standard OA method and the LP/NLP B\&B method are more efficient than their convexified counterparts, as both bound tightening and convexification cuts entail additional processing time. Nonetheless, the bound tightening and convexification cuts generally enhance the performance of both the OA and LP/NLP B\&B methods.

\begin{figure}
    \centering
    \includegraphics[width=\linewidth]{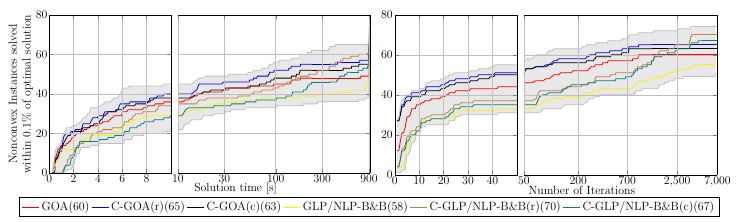}
    \caption{Time and iteration performance profile of 181 nonconvex MINLP instances}
    \label{fig:nonconvex performance profiles}
\end{figure}

The benchmark results of nonconvex instances are presented in Figure \ref{fig:nonconvex performance profiles}. Similar to the convex cases, the convexification-based OA and LP/NLP B\&B methods outperform the others in solution time and number of iterations. This consistent performance across convex and nonconvex MINLP problems demonstrates the effectiveness of the proposed convexification-based OA and LP/NLP B\&B methods.

\section{Acknowledgements}
D.B.N. was supported by the NASA Academic Mission Services, Contract No. NNA16BD14C. D.B.N. and Z.P. acknowledge the support of the startup grant of the Davidson School of Chemical Engineering at Purdue University. The authors acknowledge the support of Prof. Nick Sahnidis in providing a special version of BARON that allowed this research. The authors thank GAMS Development Copr. for providing a license to their software and various optimization solvers used for this manuscript.

\bibliographystyle{unsrtnat}
\bibliography{references}  






\end{document}